\documentclass[12pt]{amsart}
\usepackage{amssymb}
\advance\textheight by \topskip
\makeatletter
\renewcommand*\subjclass[2][2000]{%
  \def\@subjclass{#2}%
  \@ifundefined{subjclassname@#1}{%
    \ClassWarning{\@classname}{Unknown edition (#1) of Mathematics
      Subject Classification; using '2000'.}%
  }{%
    \@xp\let\@xp\subjclassname\csname subjclassname@#1\endcsname
  }%
}
\def\cal{\mathcal}
\def\Bbb{\mathbb}

\newenvironment{pf*}[1]{\proof[#1]}{\endproof}
\newcommand{\rom}{\textup}
\makeatother
\newtheorem{Theorem}[equation]{Theorem}

\newtheorem{Lemma}[equation]{Lemma}
\newtheorem{Proposition}[equation]{Proposition}

\theoremstyle{definition}

\makeatletter
\renewcommand\section{\@startsection{section}{1}%
  {\z@}{.7\linespacing\@plus\linespacing}{.5\linespacing}%
  {\reset@font\normalfont\bfseries\centering}}
\makeatother

\theoremstyle{remark}
\newtheorem{Remark}[equation]{Remark}

\newtheorem*{Acknowledgements}{Acknowledgements}
\numberwithin{equation}{section}
\numberwithin{figure}{section}

\newcommand{\thmref}[1]{Theorem~\ref{#1}}

\newcommand{\lemref}[1]{Lemma~\ref{#1}}
\newcommand{\propref}[1]{Proposition~\ref{#1}}
\newcommand{\remref}[1]{Remark~\ref{#1}}

\newcommand{\tabref}[1]{Table~\ref{#1}}

\newcommand{\Romnum}[1]{\expandafter\uppercase\expandafter{\romannumeral #1}}
\newcommand{\C}{{\Bbb C}}
\newcommand{\Z}{{\Bbb Z}}

\newcommand{\R}{{\Bbb R}}

\newcommand{\CP}{\operatorname{\C P}}

\newcommand{\rank}{\operatorname{rank}}
\newcommand{\Sign}{\operatorname{Sign}}

\newcommand{\ind}{\mathop{\text{\rm ind}}\nolimits}
\newcommand{\tr}{\operatorname{tr}}

\newcommand{\SW}{\operatorname{SW}}
\newcommand{\Spin}{\operatorname{Spin}}
\newcommand{\D}{\cal D}
\font\emailfont=cmtt10
\begin{document}
\title[Pseudofree $\Z/3$-actions on  $K3$ surfaces]
{Pseudofree $\Z/3$-actions on  $K3$ surfaces}
\author{Ximin Liu}
\address{Department of Applied Mathematics, Dalian University of Technology, Dalian 116024, China (\emailfont{xmliu@ms.u-tokyo.ac.jp})}

\author{Nobuhiro Nakamura}
\address{Research Institute for Mathematical Sciences, Kyoto university, Kyoto, 606-8502, Japan (\emailfont{nakamura@kurims.kyoto-u.ac.jp})}
\begin{abstract}
In this paper, we give a weak classification of locally linear pseudofree actions of the cyclic group of order $3$ on a $K3$ surface, and prove the existence of such an action which can not be realized as a smooth action on the standard smooth $K3$ surface.
\end{abstract}
%
%
\maketitle
%
\section{Introduction}\label{sec:intro}

Let $G$ be the cyclic group of order $3$ ($G=\Z/3$), and
suppose that $G$ acts locally linearly and pseudofreely on a $K3$ surface $X$. 
(An action on a space is called {\it pseudofree} if it is free on the compliment of a discrete subset.)  
The purpose of this paper is to give a weak classification of such $G$-actions and to prove that there exists such an action on $X$ which can not be realized by a smooth action for the standard smooth structure on $X$.
\begin{Theorem}\label{thm:main0}
There exists a locally linear pseudofree $G$-action on a $K3$ surface $X$ which can not be realized by a smooth action for the standard smooth structure on $X$.  
\end{Theorem}

After submitting this paper to the journal, the authors found that the $G$-action in \thmref{thm:main0} is unsmoothable for infinitely many smooth structures on $X$. 
This is proved in \remref{rem:inf}.

To state the result more precisely, we prepare notation. 
Let $b_i$ be the $i$-th Betti number of $X$, and $b_+$ (resp. $b_-$) be the rank of a maximal positive (resp. negative) definite subspace $H^+(X;\R)$ (resp. $H^-(X;\R)$) of $H^2(X;\R)$.
For any $G$-space $V$, let $V^G$ be the fixed point set of the $G$-action.
Let $b_\bullet^G = \dim H^\bullet (X;\R)^G$, where $\bullet = 2, +, -$.
The Euler number of $X$ is denoted by $\chi (X)$ and the signature of $X$ by $\Sign(X)$.

When we fix a generator $g$ of $G$, the representation at a fixed point can be described by a pair of nonzero integers $(a,b)$ modulo $3$ which is well-defined up to order and changing the sign of both together.
Hence, there are two types of fixed points.
\begin{itemize}
\item The type ($+$): $(1,2) = (2,1)$.
\item The type ($-$): $(1,1) = (2,2)$.
\end{itemize}
Let $m_+$ be the number of fixed points of the type ($+$), and $m_-$ be the number of fixed points of the type ($-$).

\thmref{thm:main0} immediately follows from the next theorem.

\begin{Theorem}\label{thm:main}
Let $G$ be the cyclic group of order $3$.
For locally linear pseudofree $G$-actions on a $K3$ surface $X$, we have the following{\rom :}
\begin{enumerate}
\item  Every locally linear pseudofree $G$-action on $X$ belongs to one of four types in \tabref{tab:actions}.
Furthermore, each of four types can be actually realized by a locally linear pseudofree $G$-action on $X$.
\begin{table}[h]
\caption{The classification of actions}
\label{tab:actions}
\begin{center}
\begin{tabular}{l|c|c|c|c|c|c|c}
Type & $\#X^G$ & $m_+$ & $m_-$ & $b_2^G$ & $b_+^G$ & $b_-^G$ & $\Sign(X/G)$ \\
\hline
$A_0$ & $6$ & $6$ & $0$ & $10$ & $3$ & $7$ & $-4$  \\
$A_1$ & $9$ & $3$ & $6$ & $12$ & $3$ & $9$ & $-6$ \\
$A_2$ & $12$ &  $0$ & $12$ & $14$ & $3$ & $11$ & $-8$ \\
\hline
$B$ & $3$ & $0$ & $3$ & $8$ & $1$ & $7$ & $-6$ \\
\end{tabular}
\end{center}
\end{table}
\item The type $A_1$ can not be realized by a smooth action on the standard smooth $K3$ surface.
\end{enumerate}
\end{Theorem}
\begin{Remark}\label{rem:1}
The assertion (1) in \thmref{thm:main} is an application of the remarkable result by A.~L.~Edmonds and J.~H.~Ewing \cite{EE} with Freedman's classification of simply-connected topological $4$-manifolds \cite{Freedman}.
\end{Remark}
\begin{Remark}
To prove the assertion (2), we use the mod $p$ vanishing theorem of Seiberg-Witten invariants by F.~Fang \cite{Fang}, with the fact that the Seiberg-Witten invariants for the canonical $\Spin^c$-structure of the standard smooth $K3$ surface is $\pm 1$.
\end{Remark}
\begin{Remark}
The type $A_0$, $A_1$ and $A_2$ are actions which act trivially on $H^+(X;\R)$.
\end{Remark}
\begin{Remark}\label{rem:std}
The type $A_0$ is realized by a smooth action on the Fermat quartic surface. 
(See \propref{prop:Fermat}. ) 
\end{Remark}
\begin{Remark}
We do not know whether  $A_2$ and $B$ can be realized by a smooth action for some smooth structure on a $K3$ surface, or not.
\end{Remark}
\begin{Remark}
K.~Kiyono proved the existence of unsmoothable locally linear pseudofree actions on the connected sums of $S^2\times S^2$ \cite{Kiyono}.
Although he also uses the Seiberg-Witten gauge theory, his method is different from ours.
It is interesting that he invokes the ``$G$-invariant  $10/8$-theorem'' instead of Seiberg-Witten invariants.
(A related paper is \cite{KL}.)
\end{Remark}
\begin{Acknowledgements}
The authors would like to express our gratitude to K.~Kiyono for introducing the paper \cite{EE} of A.~L.~Edmonds and J.~H.~Ewing to us.
It is also a pleasure to thank M.~Furuta for invaluable discussions.
\end{Acknowledgements}

\section{The proof of the assertion (1)}\label{sec:proof1}
As mentioned in \remref{rem:1}, the proof of the assertion (1) of \thmref{thm:main} will rely on the realization theorem by A.~L.~Edmonds and J.~H.~Ewing \cite{EE}. First, we summarize their result in the very special case when $G=\Z/3$.
\begin{Theorem}[\cite{EE}]\label{thm:EE}
Let $G$ be the cyclic group of order $3$.
Suppose that one is given a fixed point data
$$
\D = \{(a_0,b_0), (a_1,b_1), \ldots, (a_n,b_n), (a_{n+1},b_{n+1})\},
$$
where $a_i, b_i \in \Z/3\setminus\{0\}$, and a $G$-invariant symmetric unimodular form
$$
\Phi\colon V\times V\to \Z,
$$
where $V$ be a finitely generated $\Z$-free $\Z[G]$-module.
Then the data $\D$ and the form $(V,\Phi)$ are realizable by a locally linear, pseudofree, $G$-action on a closed, simply-connected, topological $4$-manifold if and only if they satisfy the following two conditions{\rom :}
\begin{enumerate}
\item The condition REP{\rom :} As a $\Z[G]$-module, $V$ splits into $F\oplus T$, where $F$ is free and $T$ is a trivial $\Z[G]$-module with $\rank_\Z T = n$.
\item The condition GSF{\rom :} The $G$-Signature Formula is satisfied{\rom :}
$$
\Sign(g, (V,\Phi)) = \sum_{i=0}^{n+1}\frac{(\zeta^{a_i} + 1)(\zeta^{b_i} + 1)}{(\zeta^{a_i} - 1)(\zeta^{b_i} - 1)},
$$
where $\zeta = \exp(2\pi\sqrt{-1}/3)$.
\end{enumerate}
\end{Theorem}
\begin{Remark}
In \cite{EE}, A.~L.~Edmonds and J.~H.~Ewing prove the realization theorem for all cyclic groups of prime order $p$, and for general  $p$, the third condition {\it TOR} which is related to the Reidemeister torsion should be satisfied.
However, when $p=3$, the condition {\it TOR} is redundant.
This follows from the fact that the class number of $\Z[\zeta]$ is $1$, and Corollary 3.2 of \cite{EE}.
\end{Remark}
Now, let us begin the proof of the assertion (1).

Suppose that a locally linear pseudofree $G$-action on $X$ is given.
First of all, the ordinary Lefschetz formula should hold: $L(g,X) = 2 + \tr ( g|_{H^2(X)}) =\#X^G$.
Noting that $\#X^G = m_+ + m_- $ and $ 2 + \tr ( g|_{H^2(X)}) \leq 24$, we obtain
$$
m_+ + m_- \leq 24.
$$
This is compatible with the condition {\it REP}.
Note that
\begin{equation*}
\chi (X/G) = \frac13 \{ 24 + 2 (m_+ + m_- )\}.
\end{equation*}

By \thmref{thm:EE}, the $G$-Signature Formula should hold:
\begin{align*}
\Sign (g,X) &= \Sign (g^2,X) = \frac13 (m_+ - m_-),\\
\Sign (X/G) &=  \frac13 \left\{ -16 + \frac23 (m_+ - m_-)\right\}.\\
\end{align*}
Since $\Sign (X/G)$ is an integer, 
$
m_+ - m_- \equiv 6 \mod 9.
$
This with the inequality $-24\leq m_+ - m_-\leq 24$ implies that
\begin{equation}\label{eq:diff}
m_+ - m_- = -21,-12,-3,6,15,24.
\end{equation}

We can calculate $b_+^G$ and $b_-^G$ from $\chi(X/G)$ and $\Sign (X/G)$.
Since $b_+^G$ is $1$ or $3$, we obtain the following:
\begin{itemize}
\item When $b_+^G=1$, $2m_+ + m_-=3$.
\item When $b_+^G=3$, $2m_+ + m_-=12$.
\end{itemize}
By these equations, \eqref{eq:diff} and non-negativity of $m_+$ and $m_-$, we obtain \tabref{tab:actions}.

Next we will prove the existence of actions. First, we construct a smooth $G$-action of type $A_0$ on the Fermat quartic surface. 
\begin{Proposition}\label{prop:Fermat}
There exists a smooth $G$-action of the type $A_0$ on the Fermat quartic surface $X$ which is defined by the equation $\sum_{i=0}^3 z_i^4 = 0$ in $\CP^3$. 
\end{Proposition}
\proof
By the symmetry of the defining equation, the symmetric group of degree $4$ acts on $X$ as permutations of variables. 
Therefore $G$ acts smoothly on $X$ via this action. 
We can easily check that the $G$-action is pseudofree, and belongs to the type $A_0$.
\endproof

To prove the existence of actions of other types, we invoke \thmref{thm:EE}. 
We need to construct $G$-actions on the intersection form.
Let $(V_{K3}, \Phi_{K3})$ be the intersection form of the $K3$ surface, which is even and indefinite.
Since an even indefinite form is completely characterized  by  its rank and signature, $(V_{K3}, \Phi_{K3})$ is isomorphic to $3H\oplus \Gamma_{16}$, where $H$ is the hyperbolic form, and $\Gamma_{16}$ is a negative definite even form of rank $16$.
We will construct $G$-actions on $3H$ and $\Gamma_{16}$ separately.

\begin{Lemma}\label{lem:E16}
For each integer $k$ which satisfies $0\leq k\leq5$, there is a $G$-action on $\Gamma_{16}$ such that
$$
\Gamma_{16}\cong (16-3k)\Z\oplus k\Z[G]\text{ as a $\Z[G]$-module}.
$$
\end{Lemma}
\proof
When $k=0$, it suffices to take the trivial $G$-action. 
Hence we suppose $k\geq 1$. 

Recall that the lattice $\Gamma_{16}$ is the set of $(x_1,\ldots,x_{16})\in(\frac12 \Z)^{16}$ which satisfy
\begin{enumerate}
\item $x_i\equiv x_j\mod \Z$ for any $i,j$,
\item $\sum_{i=1}^{16}x_i\equiv 0\mod 2\Z$.
\end{enumerate}
The unimodular bilinear form on $\Gamma_{16}$ is defined by $-\sum_{i=1}^{16}x_i^2$.

Note that the symmetric group of degree $16$ acts on $\Gamma_{16}$ as permutations of components.
For a fixed generator $g$ of $G$, define the $G$-action on $\Gamma_{16}$ by
$$g = (1,2,3)(4,5,6)\cdots(3k-2,3k-1,3k),$$
where $(l,m,n)$ is the cyclic permutation of $(x_l, x_m ,x_n)$.

As a basis for $\Gamma_{16}$, we take 
\begin{equation*}
f_i = \left\{
\begin{aligned}
e_i+e_{16},& \quad\qquad (i=1,\ldots ,9), \\
e_i-e_{16},& \quad\qquad (i=10,\ldots,15),\\
\frac12 (e_1+e_2&+\cdots+e_{16}),\, (i=16),
\end{aligned}\right.
\end{equation*}
where $e_1,\ldots,e_{16}$ is the usual orthonormal basis for $\R^{16}$. 
Then the basis $(f_1,f_2,\ldots,f_{16})$ gives required direct splitting.
\endproof
\begin{Lemma}\label{lem:3H}
There is a $G$-action on $3H$ such that $3H \cong\Z[G]\oplus\Z[G]$ as a $\Z [G]$-module, and $G$-fixed parts of a maximal positive definite subspace and a negative one of $3H\otimes\R$ both have rank $1$.
\end{Lemma}
\proof
Such a $G$-action is given as permutations of three $H$'s.
\endproof

With \lemref{lem:E16} and \lemref{lem:3H} understood, for each of $A_1$, $A_2$ and $B$, the corresponding $G$-action  on $(V_{K3}, \Phi_{K3})$ can be constructed. That is,
\begin{itemize}
\item for $A_1$, $3H \cong 6\Z$ and $\Gamma_{16}\cong \Z\oplus 5\Z[G]$, 
\item for $A_2$, $3H \cong 6\Z$ and $\Gamma_{16}\cong 4\Z\oplus 4\Z[G]$,
\item for $B$,   $3H \cong \Z [G]\oplus\Z [G]$ and $\Gamma_{16}\cong \Z\oplus 5\Z[G]$.
\end{itemize}
Now the conditions {\it REP} and {\it GSF} are satisfied. Therefore we have a locally linear pseudofree $G$-action on a closed simply-connected $4$-manifold $X$ whose intersection form is just $(V_{K3}, \Phi_{K3})$ by \thmref{thm:EE}.
Since $X$ is simply-connected and its intersection form is even, we see that $X$ is homeomorphic to the $K3$ surface by Freedman's theorem \cite{Freedman}.
Thus the assertion (1) is proved.

\begin{Remark}
By using Theorem 1.3 in \cite{BW}, we can prove that the topological conjugacy class of actions of the type $B$ is unique, that is, any action of the type $B$ is conjugate to the action which we have constructed. 
\end{Remark}
\begin{Remark}
We can also construct a locally linear pseudofree action of the type $A_0$ by \thmref{thm:EE}. 
For this purpose, we need to construct a $G$-action on $3H$ such that $3H \cong 3\Z \oplus\Z [G]$ as a $\Z[G]$-module, and the rank of a $G$-fixed maximal positive definite subspace of $3H\otimes\R$ is $3$ and the rank of a negative one is $1$. 
Such a $G$-action on $3H$ is constructed from the cohomology ring of a $4$-torus with a $G$-action as follows:

Let $\zeta =\exp(2\pi\sqrt{-1}/3)$, and consider the lattice $\Z\oplus\zeta\Z\subset\C$.
For each $i = 0,1,2$, let us consider a $2$-torus $T_{\zeta^i} = \C/(\Z\oplus\zeta\Z)$ with a $G$-action, where the $G$-action is defined by the multiplication by $\zeta^i$.
Next, consider the $4$-torus $T_{12}=T_{\zeta}\times T_{\zeta^2}$ with the diagonal $G$-action.
Then we can prove that the induced $G$-action on $H^2(T_{12};\Z)$ has required properties.%

Using this with a $G$-action on $\Gamma_{16}$ such that $\Gamma_{16}\cong \Z\oplus 5\Z[G]$, we obtain a $G$-action of the type $A_0$ by \thmref{thm:EE}.  
\end{Remark}

\section{The proof of the assertion (2)}\label{sec:proof2}

In this section, we consider $X$ as the smooth $K3$ surface with the standard smooth structure.
Suppose now that a smooth action of the type $A_1$ exists.
To obtain a contradiction, we use a Seiberg-Witten invariant of $X$.
Recall that, for a smooth $4$-manifold with $b_1=0$ and $b_+\geq 2$, Seiberg-Witten invariants constitute a map from the set of equivalence classes of $\Spin^c$-structures on $X$ to $\Z$.
That is, for a $\Spin^c$-structure $c$, the corresponding Seiberg-Witten invariant $\SW_X(c)$ is given as an integer.

We use the canonical $\Spin^c$-structure $c_0$ which is characterized as one whose determinant line bundle $L$ is trivial in the case of $K3$ surface $X$.
Note that $c_0$ is also characterized as the $\Spin^c$-structure which is determined by the $\Spin$-structure.

Since $X$ is simply-connected and $L$ is trivial, we can see that every $G=\Z/3$-action on $X$ lifts to a $G$-action on the $\Spin^c$-structure $c_0$.
Then, the $G$-index of the Dirac operator $D_X$ can be written as $\ind_G D_X = \sum_{j=0}^2 k_j\C_j \in R(G) \cong \Z[t]/(t^3=1)$, where $\C_j$ is the complex $1$-dimensional weight $j$ representation of $G$ and $R(G)$ is the representation ring of $G$.

F.~Fang \cite{Fang} proves the mod $p$ vanishing theorem under a $\Z/p$-action where $p$ is a prime.
\begin{Theorem}[\cite{Fang}]\label{thm:Fang}
Let $Y$ be a smooth closed oriented $4$-dimensional $\Z/p$-manifold with $b_1=0$ and $b_+\geq 2$, where $p$ is a prime. Suppose that $c$ is a $\Spin^c$-structure on which the $\Z/p$-action lifts, and that $\Z/p$ acts trivially on $H^+(Y;\R)$.
If $2k_j \leq b_+ -1$ for $j=0,\ldots,p-1$, then
\begin{equation*}
\SW_Y(c) \equiv 0 \mod p.
\end{equation*}
\end{Theorem}
\begin{Remark}
The second author generalized \thmref{thm:Fang} to the case when $b_1>0$ \cite{Nakamura}.
\end{Remark}

On the other hand, it is well-known that $\SW_X(c_0) = \pm 1$ for the standard $K3$ surface $X$. (See e.g. \cite{FM} or \cite{T}.)
Therefore, in the case when $G$ acts on $(X,c_0)$, we have $k_j >1$ for some $j$ by \thmref{thm:Fang}.

Coefficients $k_j$ are calculated by the $G$-spin theorem. 
(For the $G$-spin theorem, we refer \cite{AB,AH,LM,Sh}.)
For the fixed generator $g\in G$, the Lefschetz number $\ind_g D_{X}$ is calculated by the formula as 
\begin{equation*}
\ind_g D_{X} =\sum_{j=0}^{2} \zeta^j k_j = \sum_{P\in X^G} \nu(P), 
\end{equation*}
where $\zeta=\exp (2\pi\sqrt{-1}/3)$ and $\nu(P)$ is a complex number associated to each fixed point $P$ given as follows.

Suppose that a fixed point $P$ has the representation type $(a, b)$ with respect to $g$. 
Then the number $\nu (P)$ associated to $P$ is given by,
\begin{equation}\label{eq:nup}
\nu (P) = \frac1{{(\zeta^{a})}^{1/2} - {(\zeta^{a})}^{-1/2}}\frac1{{(\zeta^{b})}^{1/2} - {(\zeta^{b})}^{-1/2}}.
\end{equation}
The signs of ${(\zeta^{a})}^{1/2}$ and ${(\zeta^{b})}^{1/2}$ are determined such that
$$
\left\{{(\zeta^{a})}^{1/2}\right\}^3 =\left\{{(\zeta^{b})}^{1/2}\right\}^3 = 1.  
$$
(This is because, in our case, the $g$-action on the $\Spin$-structure generates a $G$-action on the $\Spin$-structure. See \cite[p.20]{AH} or \cite[p.175]{Sh}.)

With the above understood, we obtain
\begin{align*}
\ind_g D_X & =k_0 + \zeta k_1 + \zeta^2 k_2 = \frac13 (m_+ - m_-),\\
\ind_{g^2} D_X & = k_0 + \zeta^2 k_1 + \zeta k_2 = \frac13 (m_+ - m_-),\\
\ind_1 D_X & = k_0 +  k_1 +  k_2 = 2.\\
\end{align*}
Solving these equations, we have
\begin{align*}
 k_0 &= \frac19 \left\{ 6 + 2 (m_+ - m_-)\right\},\\
k_1 = k_2 &=  \frac19 \left\{ 6 - (m_+ - m_-)\right\}.
\end{align*}
In the case of an action of type $A_1$, $m_+=3$ and $m_- =6$.
Hence, we have $k_0=0$ and $k_1=k_2=1$.
Therefore there is no $j$ so that $k_j >1$. This is a contradiction. Thus the assertion (2) is proved.
\begin{Remark}\label{rem:inf}
It is clear that a proposition similar to (2) of \thmref{thm:main} is true for the smooth structure such that the Seiberg-Witten invariant for the $\Spin^c$-structure with trivial determinant line bundle is not congruent to $0$ modulo $3$. 
Let us examine elliptic surfaces which are homeomorphic to $K3$. 
Consider relatively minimal regular elliptic surfaces with at most two multiple fibers whose Euler number is $24$. 
Let $p$ and $q$ be the multiplicities of multiple fibers, and let us write such elliptic surface as $E(2)_{p,q}$. 
(We assume that $p$ and $q$ may be $1$.)
The following are known about $E(2)_{p,q}$.
\begin{enumerate}
\item $E(2)_{1,1}$ (no multiple fiber) is diffeomorphic to the standard $K3$.
\item $E(2)_{p,q}$ is homeomorphic to the $K3$ surface if and only if $\gcd(p,q)=1$. (See e.g.\cite{Ue}.)
\item $E(2)_{p,q}$ is not diffeomorphic to $E(2)_{p^\prime,q^\prime}$ if $pq\neq p^\prime q^\prime$\cite{FM0}. 
\item Let $c_0$ be the $\Spin^c$-structure with trivial determinant line bundle. 
If $p$ and $q$ are odd, then $\SW_{E(2)_{p,q}}(c_0)=\pm 1$ \cite{FM2,FS}. 
\end{enumerate}
Thus we see that the type $A_1$ can not be realized by a smooth action on  $E(2)_{p,q}$ such that $\gcd(p,q)=1$ and $p$ and $q$ are odd. 
Note that there are infinitely many $(p,q)$ which give different smooth structures.
\end{Remark}

\end{document}